\documentclass[review]{elsarticle}

\usepackage{hyperref}
\usepackage{epsfig}
\usepackage{graphicx}
\usepackage{xcolor}
\usepackage{amsmath,amssymb}
\newcommand{\R}{{\mathbb{R}}}
\newcommand{\N}{{\mathbb{N}}}
\newcommand{\Z}{{\mathbb{Z}}}
\newcommand{\C}{{\mathbb{C}}}
\newcommand{\supp}{\operatorname{supp}}
\newcommand{\card}{\operatorname{card}}

\newcommand{\zb}[1]{\ensuremath{\boldsymbol{#1}}}
\newcommand{\e}{\mathrm{e}}
\newcommand{\im}{\mathrm{i}}
\newcommand{\dx}{\mathrm{d}}
\newcommand{\spa}{\mathrm{span}}
\newcommand*{\vv}[1]{\vec{\mkern0mu#1}}
\newcommand{\vc}[1]{{\overrightarrow{#1}}}

\begin{document}

\begin{frontmatter}

\title{A direct solver for the phase retrieval problem in ptychographic imaging}


\author[TUM,HMGU]{Nada Sissouno\corref{mycorrespondingauthor}}
\cortext[mycorrespondingauthor]{Corresponding author}
\ead{sissouno@ma.tum.de}
\author[HIT]{Florian Bo\ss mann}
\author[HMGU]{Frank Filbir}
\author[MSU]{Mark Iwen}
\author[DESY]{Maik Kahnt}
\author[UCSD]{Rayan Saab}
\author[DESY]{Christian Schroer}
\author[HMGU]{Wolfgang zu Castell}

\address[TUM]{Technical University of Munich, Faculty of Mathematics, Boltzmannstra\ss e 3, Garching, 85748, Germany}
\address[HMGU]{Helmholtz Zentrum M\"unchen, Scientific Computing Research Unit, Neuherberg,
85674, Germany}
\address[DESY]{Deutsches Elektronen Synchrotron DESY, Hamburg, 22607, Germany}
\address[HIT]{Harbin Institute of Technology, Department of Mathematics, Harbin 150001, China}
\address[MSU]{Michigan State University, Department of Mathematics and Department of CMSE,
East Lansing, MI, 48824, USA}
\address[UCSD]{University of California San Diego, Department of Mathematics, La Jolla, CA,
92093, USA}
%%%%%%%%%%%%%%%%%%%%%%%%%%%%%%%%%%%%%%%%%%%%%%
\begin{abstract}
Measurements achieved with ptychographic imaging are a special
case of diffraction measurements. They are generated by
illuminating small parts of a sample with, e.g., a focused X-ray beam.  By shifting
the sample, a set of far-field diffraction patterns of the
whole sample are then obtained. From a mathematical point of view those
measurements are the squared modulus of the windowed Fourier transform of the sample.
Thus, we have a phase retrieval problem for local Fourier measurements.
A direct solver for this problem was introduced by Iwen, Viswanathan and Wang in
2016 and improved by Iwen, Preskitt, Saab and Viswanathan in 2018.
Motivated by the applied perspective of ptychographic imaging,
we present a generalization of this method and compare the different versions
in numerical experiments. The new method proposed herein turns out to be more stable,
particularly in the case of missing data.
\end{abstract}
%%%%%%%%%%%%%%%%%%%%%%%%%%%%%%%%%%%%%%%%%%%%%%

\end{frontmatter}

%\linenumbers
%%%%%%%%%%%%%%%%%%%%%%%%%%%%%%%%%%%%%%%%%%%%%%
\section{Introduction}
%%%%%%%%%%%%%%%%%%%%%%%%%%%%%%%%%%%%%%%%%%%%%%
Ptychography refers to a diffraction imaging technique which collects a number of diffraction patterns of an object in the far-field, where each pattern is 
generated by illuminating a small subregion one at a time \cite{HeHo70, Ro08, SeEtal17, Pe18}. The selection of the small and necessarily overlapping subregions is managed by using a mask, or window, 
placed between the X-ray beam and the object. For every shift of the mask, a diffraction pattern is recorded and the imaging task consists of reconstructing 
the object function from this collection of measurements.
Since the measurement takes place in the far-field it is given 
as the squared modulus of the windowed Fourier transform of the object. This, in particular, results in a loss of the phase information of the signal. For the 
image formation we therefore  have to reconstruct the object function from phaseless localized Fourier transform data. The imaging task in ptychography can thus be formulated in 
mathematical terms as follows. For an unknown object modeled as a complex-valued function $f\in L^2(\R^d)$ supported on a compact set $\Omega\subset\R^d$  
we measure 
\begin{equation}\label{eq:1.1}
y(\tau,\omega)=\big|\mathcal{F}[f\, T_\tau w](\omega)\big|^2 =\big|\mathcal{F}[T^\ast_\tau f\, w](\omega)\big|^2,
\end{equation}
where $w\in L^2(\R^d)$ is a compactly supported window function, $T_\tau f(x)=f(x-\tau)$ is the translation operator with adjoint 
$T^\ast_\tau f(x)=f(x+\tau)$, and $\mathcal{F}$ is the Fourier transform defined as $\mathcal{F}f(\omega)=\int_{\R^d}f(x)\, \e^{-2\pi\im\,  \omega\cdot x}\, \dx x$. 
We have to invert the non-linear mapping $f\to y$. Note that the latter formulation in \eqref{eq:1.1} more closely represents the real experimental situation where the object is 
shifted instead of the window. Of course, in concrete applications we are given only samples $\{y(\tau_\ell,\xi_k)\}$ and the inversion problem 
has to be formulated in a fully discrete regime. Discretization of \eqref{eq:1.1} on a grid $\{n/N:n=0,\dots, N-1\}$ leads to 
\begin{equation}\label{eq:1.2}
y_{\ell,k}=|\langle \zb f,M_kT_\ell \zb w\rangle|^2, 
\end{equation}
where $\zb f=(f(0),\dots, f(N-1))^T$ and $\zb w=(w(0),\dots, w(N-1))^T\in\C^N$ are vectors containing the sampled values of the function $f$ resp. $w$. 
The translation operator $T_\ell$ and modulation operator $M_k$ are acting on the entries of the vectors by 
$T_\ell f(n)=f(n-\ell\mod N)$ resp. $M_kf(n)=\e^{2\pi\im nk/N} f(n)$. They have obvious matrix representations.
Recovery of $\zb f$ from data \eqref{eq:1.2} falls in the class 
of problems where a vector $\zb x$ has to be recovered from data of the form $|\langle \zb x,\zb\varphi_m\rangle|$ with $\{\zb \varphi_m\}$ being a frame for 
$\C^N$. We are considering here the case where frame is given as a discrete Gabor frame $\{M_kT_\ell\zb w\}$ for $k=0,\dots,K$ and {$\ell\in\mathcal{L}\subseteq \{0,\dots,N-1\}$ with $\card(\mathcal{L})=L+1$ for} $K,\,L\leq N-1$. 

The general Phase Retrieval Problem has a long history and it was tackled in many different ways and under various assumptions. 
We will  make no attempt to review these developments in detail here but refer to \cite{Lu17, ElHaMi16} and references cited there. Among the diverse 
techniques one method seems to be 
particularly popular among practitioners. In the physics literature this method is known as the Ptychographic Interative Engine (PIE), and is regarded as a standard approach for image formation from ptychographic data within the physics community \cite{MaRo09, MaJoLi17}. 
Mathematically this approach is an alternating projection method which goes back to the work of Gerchberg and Saxton \cite{GeSa72} and Fienup 
\cite{Fi86}. The method projects alternately on the set which consists of functions with support in $\Omega$ and the set of functions which agree with the 
measurement \cite{BaCoLu02}. The latter set is non-convex which makes the problem notoriously difficult
to analyze theoretically. Due to the non-convexity, the alternating projection method can {converge to a stationary point that differs from the true solution}.  Determining a good starting point in general is also not easy and some attempts were made to come up with a 
good initial guess, see for example \cite{MaTuWu16} for a recent approach in this regard.  

The method we present in this paper is an adoption of a fast direct solver for the phase retrieval problem \eqref{eq:1.2} as it was developed by Iwen et. al. in \cite{IwViWa16,IwPrSaVi18}. 
The method is based on a lifting scheme as used in the PhaseLift algorithm \cite{CaStVo13} which transforms the discretized non-linear problem \eqref{eq:1.2} into a linear 
problem for the lifted variables $\zb f\zb f^\ast$. After recovering some of the entries of the lifted variables $\zb f\zb f^\ast$ the phase of the individual entries of $\zb f$ can be 
determined by an angular synchronization approach. This finally results in the reconstruction of the function $f$ on the grid up to a global phase multiple provided that the measurements are sufficiently informative. 

In this paper we will demonstrate how to apply a new algorithm for phase retrieval from short-time Fourier measurements to the concrete experimental setup of
ptychography.  The paper is organized as follows. We state the Algorithm in the Section~\ref{sec:2}. Here, we first introduce the 1D case and later give its 2D version. In Section~\ref{sec:Ex} we demonstrate the method with numerical examples. Finally, Section~\ref{sec:Future} concludes with a discussion of the proposed technique.
%%%%%%%%%%%%%%%%%%%%%%%%%%%%%%%%%%%%%%%%%%%%%%
\section{Phase Retrieval from Localized Fourier Measurements}\label{sec:2}
%%%%%%%%%%%%%%%%%%%%%%%%%%%%%%%%%%%%%%%%%%%%%%
\subsection{Description of the Algorithm}\label{ssec:2.1}
%%%%%%%%%%%%%%%%%%%%%%%%%%%%%%%%%%%%%%%%%%%%%%
Let $f,w\in L^2(\R^d)$ be compactly supported functions. Without loss of generality we may assume $\supp(f)\subseteq [0,1]^d$. The short-time Fourier transform 
of $f$ with window $w$ is defined as 
\begin{equation}\label{eq:2.1}
V_w f(\tau,\omega):=\int_{\R^d} f(t)\, \bar{w}(t-\tau)\, \e^{-2\pi\im \omega\tau}\, \dx t=\langle f,M_\omega T_\tau w\rangle_{L^2(\R^n)},
\end{equation}
where $T_\tau w(t)=w(t-\tau),\ M_\omega w(t)=\e^{2\pi\im \omega t}w(t)$ are the translation resp. modulation operator. We concentrate on the cases $d\in\{1,2\}$. For clarity of presentation, we first 
restrict ourselves to the case $d=1$. The case $d=2$ will be addressed in Subsection~\ref{ssec:2.2}.  For discretization let $N\in\N$ and consider the grid
$\Gamma=\{n/N:n=0,\dots N-1\}$. Discretization of the Fourier integral in \eqref{eq:2.1} on $\Gamma$ and subsequential evaluation of the resulting semi-discrete 
transform w.r.t. $\tau$ on the grid $\Gamma$ leads to the fully discretized transform which can be considered as a short-time Fourier transform on the cyclic group 
$\Z_N$. It is given as 
\begin{equation}\label{eq:2.2}
V_w f(\ell,k)=\frac{1}{N}\sum_{n\in\Z_N}f(n) \bar{w}(m-\ell)\, \e^{-2\pi\im k\cdot n/N},\quad \ell,k\in\Z_N.
\end{equation}
This can be expressed in a more condensed form as 
\begin{equation}\label{eq:2.3}
V_w f(\ell,k)=\frac{1}{N}\langle\, \zb f,M_kT_\ell\zb w\rangle_{\C^N}.
\end{equation}
where $\zb f,\,\zb w\in \mathbb{C}^N$,
$T_\ell f(n)=f(n-\ell\mod N)$, and $M_kf(n)=\e^{2\pi\im nk/N} f(n)$ as already defined after
\eqref{eq:1.2}. The scaling factor $1/N$ in \eqref{eq:2.3} is of no relevance for our consideration and will 
therefore be neglected henceforth. \\
The reconstruction problem can now be formulated as follows. We have to reconstruct the vector $\zb f$ from data 
\begin{equation}\label{eq:2.4}
y_{\ell,k}=\big|\langle \zb f,T_\ell\bar{\zb w}_k\rangle_{\C^N}\big|^2 ,\quad k=0,\dots, K,\ {\ell\in\mathcal{L}\subseteq\{0,\dots,N-1\}}, 
\end{equation}
{with $\card(\mathcal{L})=L+1$ and} $K,L\leq N-1$, where $\zb w_k=M_k\zb w$ is the modulation of the vector of window $\zb w$. We will henceforth assume that $\zb w$ is supported 
on the first $s$ entries, i.e., $w(n)\not= 0$ for $0\leq n< s$ and $0$ elsewhere. \\ 
These nonlinear measurements can be ``lifted'' to linear measurements on the space of matrices as follows 
\begin{equation}\label{eq:2.5}
\big|\langle \zb f,T_\ell\bar{\zb w}_k\rangle_{\C^N}\big|^2=\langle \zb f,T_\ell\bar{\zb w}_k\rangle_{\C^N}\ \overline{\langle \zb f,T_\ell\bar{\zb w}_k\rangle_{\C^N}}= 
\langle\zb f\zb f^\ast,T_\ell\bar{\zb w}_k(T_\ell\bar{\zb w}_k)^\ast\rangle_{HS}.
\end{equation}
Here $\langle\cdot,\cdot\rangle_{HS}$ is the Hilbert-Schmidt inner product $\langle \zb A,\zb B\rangle_{HS}=\mathrm{trace}(\zb A^\ast \zb B)$.\\
For $\zb X\in\C^{N\times N}$ we arrange the numbers $\langle \zb X,T_\ell\bar{\zb w}_k(T_\ell\bar{\zb w}_k)^\ast\rangle_{HS}$ 
as a vector in $\C^D$ for $D=(K+1)\, (L+1)$ and define a linear operator
$\mathcal{A}:\C^{N\times N}\to \C^D$ by
\begin{equation}\label{eq:2.6}
\mathcal{A}(\zb X)=\big(\langle \zb X,T_\ell\bar{\zb w}_k(T_\ell\bar{\zb w}_k)^\ast\rangle_{HS}\big)_{\alpha=1}^D,\quad 
\alpha=k+1+\ell(K+1)
\end{equation}
and we will call this operator the measurement operator as it coincides with the vector of measurements if $\zb X=\zb f\zb f^\ast$ according to \eqref{eq:2.5}. The 
operator \eqref{eq:2.6} can not be injective on $\C^{N\times N}$ when $\zb w$ is supported on the first $s$ entries of $N$. But,
depending on the choice of $\zb w$, it might be stably invertible if we restrict it to the space 
\begin{equation}\label{eq:2.7}
\spa\{T_\ell\bar{\zb w}_k\bar{\zb w}_k^\ast T_\ell^\ast:
{\ell\in\mathcal{L}},\ k=0,\dots, K\}=:\mathbb{P}_{LK}.
\end{equation}
We denote the corresponding projection operator by $\mathcal{T}_{LK}:\C^{N\times N}\to \mathbb{P}_{LK}$. We will
refer to this operator as the {\em tight projector}.
Let $\mathcal{A}_{LK}=\mathcal{A}|_{\mathbb{P}_{LK}}$ be the restriction of $\mathcal{A}$ to the space $\mathbb{P}_{LK}$.  
Depending of the choice
of $\zb w$, the matrix of this linear operator with respect to the generating system $\{T_\ell\bar{\zb w}_k\bar{\zb w}_k^\ast T_\ell^\ast\}_{\ell,k}$ is a symmetric positive definite matrix. Hence, the operator is (stably) invertible and we can determine $\zb X=\mathcal{T}_{LK}(\zb f\zb f^\ast)$ from the data uniquely solving 
\begin{equation}\label{eq:2.9}
\mathcal{A}_{LK}(\zb X)=\zb y,\quad \zb y=(y_{\ell,k})
\end{equation}
since
\begin{equation}\label{eq:2.8} 
y_{\ell,k}=\langle\mathcal{T}_{LK}(\zb f\zb f^\ast), T_\ell\bar{\zb w}_k\bar{\zb w}_k^\ast T_\ell^\ast\rangle_{HS}.
\end{equation}

Note that for $L=N-1$ we have 
\begin{equation*}%\label{eq:2.7}
\mathbb{P}_{LK}=\spa\{\zb E_{i,j}:|i-j\mod N|<s\}=:\mathbb{V}_s,
\end{equation*}
where the $\zb E_{i,j}$'s constitute the standard basis of $\C^{N\times N}$. This case coincides with the projection considered in \cite{IwPrSaVi18}. 
We have $\dim\mathbb{V}_s=(2s-1)N$. Thus, invertibility can only be achieved in general for $K+1=2s-1$. {If $L<N-1$, $\mathbb{P}_{LK}$ can not be expressed in terms
of $(L+1)(K+1)$ standard basis vectors $\zb E_{i,j}$.} 
Therefore, the algorithm in \citep{IwPrSaVi18}
needs to be adapted in our context.
In Section~\ref{sec:Ex} we will compare the two methods where we will refer to
the projection operator given in terms of $\zb E_{i,j}$'s in \eqref{eq:2.9} and \eqref{eq:2.8}
as {\em pattern projector} $\mathcal{P}_{LK}$, since the projector
reflects the zero pattern generated by the support of the windows. 

In order to determine an approximation of
the vector $\zb f$ from $\zb X=\mathcal{T}_{LK}(\zb f\zb f^\ast)$
we consider the representation $\zb Z$ of $\zb X$ in the standard basis. 
The amplitudes can now be determined by taking the
square-roots of the main-diagonal of $\zb Z$. To
reconstruct the phases we can use the following angular synchronization technique: Let $\zb z\in\C^N$ with $z_n=|z_n|\e^{\im\Theta_n}$, 
let $\tilde{\zb z}:=\zb z/|\zb z|=(\e^{\im\Theta_n})_n$ and define $\tilde{\zb Z}:=\zb Z/|\zb Z|$, where the operations are
considered elementwise for the non-zero entries. The phases
of $\zb z$ are then given by the first eigenvector of $\tilde{\zb Z}$. 
In case of $L=N-1$ it can be shown \cite{ViIw15} that 
the reconstruction is exact. In other words, up to a global phase shift, the phases of $\zb z$ coincide with the phases of $\zb f$.
For $L<N-1$ this process can be stabilized by considering the first eigenvector of $\zb D^{-1/2}\tilde{\zb Z}\zb D^{-1/2}$. Here $\zb D$ is the degree matrix which is a diagonal matrix whose diagonal entries coincide with the number of non-zero entries in the corresponding row of $\tilde{\zb Z}$ (see \cite{Pr18} for more details).\\

The algorithm for recovering $\zb f$ from measurements $y_{\ell,k}$ consists of five steps which can be summarized as follows.

\noindent \begin{minipage}{\textwidth}
\noindent\rule{\textwidth}{0.2pt}

{\bf Algorithm}\vspace*{-6pt}

\noindent\rule{\textwidth}{0.2pt}\vspace*{10pt}

 \begin{tabular}{ll}
 {\bf Input:}& Measurements $\zb y=(y_{\ell,k})\in\C^D$\\[0.5ex]
 {\bf Output:}& $\tilde{\zb f}\approx \zb f$\\[0,5ex]
 1.& Compute $\zb X=\mathcal{A}_{LK}^{-1}(\zb y)$.\\
 2.& Calculate $\zb Z$, the representation of $\zb X$ in the standard basis.\\
 3.& Form $\tilde{\zb Z}=\zb Z/|\zb Z|$, normalizing non-zero entries of $\zb Z$.\\
 4.& Compute the eigenvector $\tilde{\zb z}$ of $\tilde{\zb Z}$ corresponding to the largest \\
  &  eigenvalue.\\
 5. & Form $\tilde{\zb f}$ via $f_n=\sqrt{Z_{n,n}}\ \tilde{z}_n$.
 \end{tabular}
 
\noindent\rule{\textwidth}{0.2pt}
\end{minipage}
\\

The crucial step of the algorithm is the first, i.e., we have to ensure by a suitable choice of the window $\zb w$ that \eqref{eq:2.9} has a unique solution 
which can be computed in a stable manner. In \cite{IwViWa16, IwPrSaVi18} it has been shown that this is the case for the following choice for $\zb w$: 
\begin{equation}\label{eq:2.11}
w_a(n)=
\left\{
             \begin{array}{cl}
             (2s-1)^{-1/4}\ \e^{-n/a},& n<s,\\[1ex]
             0,& otherwise,
             \end{array}
\right.             
\end{equation}
where $a\in[4,\infty)$. In the numerical examples we will call these {\em exponential windows} (EW).

However, this window does not reflect the concrete experimental situation. Closer to physical reality is a shape of the window function which is given by an  
Airy function, i.e., $w(t)=(J_1(t)/t)^2$ where $J_1$ is the Bessel function of first order. Experiments show that a Gaussian function seems to be an 
acceptable approximation to this type of window function (see, e.g., \cite{SeEtal17}).
Further, the Gaussian function should be $(i)$ centered at the midpoint of the support and $(ii)$
normalized such that the norm coincides the 
number of photons in the experiment denoted by $n(p)$. 
Therefore, we propose to use the {\em Gaussian windows} (GW) constructed as follows: we consider a Gaussian window function
with normalization $c_{n(p)}$ such that $(ii)$ is satisfied.
To ensure $(i)$, we choose the support according to an 
$\alpha$-quantile $t_\alpha$ for some $\alpha\in (0,1)$, i.e., the continuous window function is given by
$w_\alpha(t)=c_{n(p)}\,e^{-{t^2}/{2}}\cdot\chi([-t_\alpha,t_\alpha])$. By uniform sampling on
$[-t_\alpha,t_\alpha]$ we get $\zb w$ with entries
\begin{align}\label{eq:GW_coef}
w_\alpha(n)=\begin{cases}
c_{n(p)}
\,e^{-\frac{t_\alpha^2}{2(s-1)^2}
(2n-s+1)^2}
& \text{for }n<s\\
0 &\text{for }n\ge s
\end{cases}.
\end{align}

In Section~\ref{sec:Ex} %\ref{sec:Ex}
we will compare reconstructions using the exponential window (EW) and Gaussian window (GW) as well as
the use of the pattern projector and tight projector for
2D data. Thus, we first give some remarks on the 2D case.
%%%%%%%%%%%%%%%%%%%%%%%%%%%%%%%%%%%%%%%%%%%%%%
\subsection{The 2D Case}\label{ssec:2.2}
%%%%%%%%%%%%%%%%%%%%%%%%%%%%%%%%%%%%%%%%%%%%%%
The discretization of \eqref{eq:2.1} for the case $d=2$ uses the cartesian grid $\Gamma\times\Gamma$ and finally leads in an analogous manner as for the 1D 
case to 
\begin{equation}\label{eq:2.12}
V_wf(\ell,k)=\frac{1}{N^2}\sum_{n_1,n_2\in\Z_N} f(n_1,n_2)\, w(n_1-\ell_1,n_2-\ell_2)\, \e^{-2\pi\im (k_1n_1+k_2n_2)/N^2}, 
\end{equation}
with $\ell=(\ell_1,\ell_2),k=(k_1,k_2)\in\Z_N\times\Z_N$. Analogously to \eqref{eq:2.3} we can express this as 
\begin{equation}\label{eq:2.13}
V_wf(\ell,k)=\frac{1}{N^2}\langle\zb F,M_kT_\ell\zb W\rangle_{HS},
\end{equation}
where $\zb F=\big(f(n_1,n_2)\big)_{n_1,n_2=0}^{N-1},\ \zb W=\big(w(n_1,n_2)\big)_{n_1,n_2=0}^{N-1}\in \C^{N\times N}$ and with the usual convention that 
the translation  resp. modulation operators are acting entrywise. As in the 1D case we will ignore the factor $1/N^2$ henceforth. \\ 
Following \cite{IwPrSaVi17}, for the 2D case we make the assumption that the window function $w$ in \eqref{eq:2.1} separates w.r.t. the variables, i.e.,  $w(t)=w(t_1,t_2)=u(t_1)\, v(t_2)$. 
We assume moreover that $\supp(u)=\supp(v)$.  Clearly,  the Gaussian window $w(t)=\e^{-|t|^2}$ satisfies this assumption, but note that the Airy function does not. With this assumption we have 
$\zb W=\zb u\zb v^\ast$ with $\zb u,\zb v\in\C^N$ and, moreover, $M_kT_\ell\zb W=(M_{k_1}T_{\ell_1}\zb u)\, (M_{k_2}T_{\ell_2}\zb v)^\ast$. In order to recover 
$\zb F$ from measurements 
\begin{equation}\label{eq:2.14}
Y_{\ell,k}=\big| \langle\zb F,M_kT_\ell\zb W\rangle_{HS}\big|^2=\big| \langle\zb F,T_\ell \zb W_k \rangle_{HS}\big|^2
\end{equation}
with $\zb W_k=M_k\zb W$ we adapt the algorithm for the 1D case as follows.
For $\zb X\in\C^{N\times N}$ let $\vv{\zb X}=\mathrm{vec}(\zb X)$ be its vectorization. Note that for $\zb X,\zb Y\in\C^{N\times N}$ we have 
$\langle \zb X,\zb Y\rangle_{HS}=\langle\vv{\zb X},\vv{\zb Y}\rangle_{\C^{N^2}}$. With this preparation it is now obvious how to transfer the algorithm presented 
in Section~\ref{ssec:2.1}
to the 2D case. Applying the same lifting step as in \eqref{eq:2.5}  we obtain 
\begin{equation}\label{eq:2.15}
Y_{\ell,k}=\big\langle\vv{\zb F}\vv{\zb F}^\ast, \vc{(T_{\ell}\zb W_k})\, \vc{(T_{\ell}\zb W_k})^\ast \big\rangle_{HS}
\end{equation}
where according to $\zb W=\zb u\zb v^\ast$, $T_\ell\zb W_k =(T_{\ell_1}M_{k_1}\zb u)\ (T_{\ell_2}M_{-k_2}\zb v)^\ast$.  

%%%%%%%%%%%%%%%%%%%%%%%%%%%%%%%%%%%%%%%%%%%%%%%%%%%%%%%%%%%%%%%%%%%%%%%%%%%%%%%%%%%%%%%%%%%
\section{Numerical examples}\label{sec:Ex}
\begin{figure}[ht]
	\center{
		\begin{tabular}{c}
			\includegraphics[height=3.2cm]{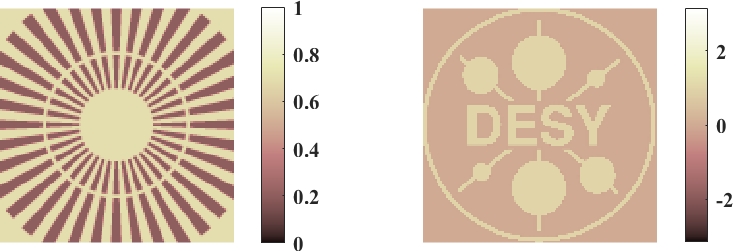}
	\end{tabular}}
	\caption{Original amplitude (left) and original phase (right) of the simulated object}
	\label{fig:originals}
\end{figure}

%%%%%%%%%%%%%%%%%%%%%%%%%%%%%%%%%%%%%%%%%%%%%%%%%%%%%%%%%%%%%%%%%%%%%%%%%%%%%%%%%%%%%%%%%%%
We test the proposed method on simulated ptychographic data using the object shown in Figure \ref{fig:originals}. {For an easier comparison the color range of all plots has been set to $[0,1]$ for amplitudes and $[-\pi,\pi]$ for phases. However, keep in mind that the phase is only reconstructed up to a global shift. For the experiments we use an} object size of $128\times128$ pixels and each pixel is about $30\times30$ nm. 
The size of the Fourier measurements per shift was chosen to be $15\times15$ pixels, i.e.,
$K=15$ frequencies in both directions. This coincides to a window size of 8x8 pixels ($s=8$ in each direction). To generate the measurement data a Gaussian beam with main-focus over the support of the window was simulated.

\begin{figure}[ht]
	\center{
		\begin{tabular}{cc}
			\includegraphics[height=3.2cm]{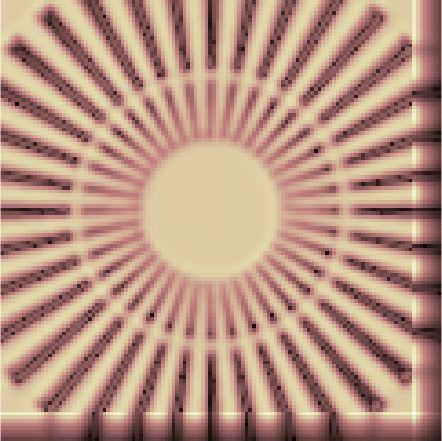} &
			\includegraphics[height=3.2cm]{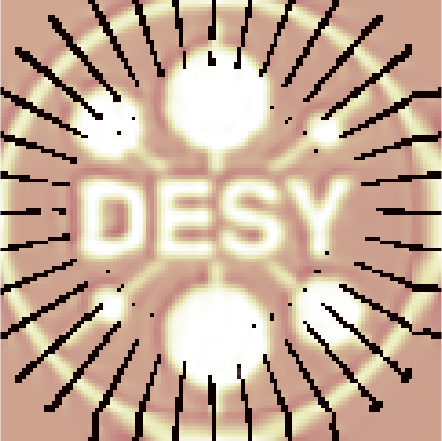} \\
			\includegraphics[height=3.2cm]{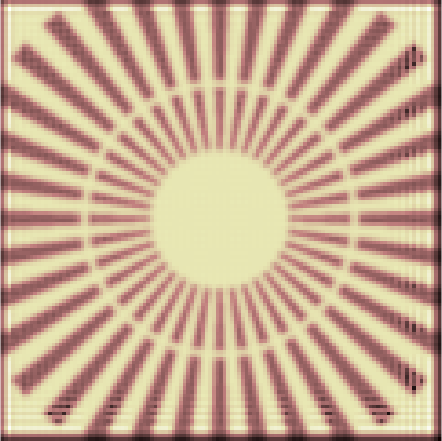} &
			\includegraphics[height=3.2cm]{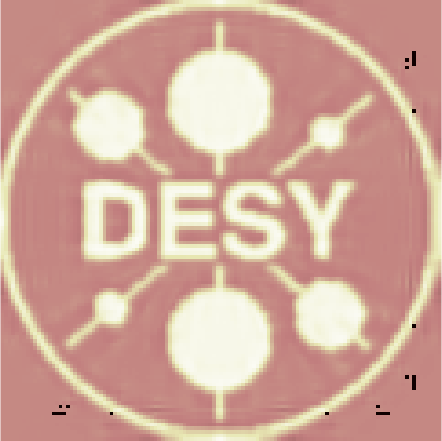}
	\end{tabular}}
	\caption{Reconstructed amplitude (left) and phase (right) using exponential window (top)  and Gaussian window (bottom)}
	\label{fig:allShifts}
\end{figure}

In a first experiment we reconstruct the object using all shifts of the window function that fit into the 128x128 pattern of the object. Given a window size of 8x8 pixels, this corresponds to $L=121$ shifts for each dimension. Note that this already differs from the setup given in \cite{IwPrSaVi18} where also circulant shifts are considered, i.e., $L=N-1=127$. We compare the reconstruction quality using the exponential window {\eqref{eq:2.11}} analyzed in \citep{IwViWa16,IwPrSaVi18} against the Gaussian window {\eqref{eq:GW_coef}} with $\alpha=0.99$, which more closely approximates experimental setups. The results are shown in Figure \ref{fig:allShifts}. Since we do not consider circulant shifts, both reconstructions show artifacts at the sides of the images. Because the exponential window is not centered, these artifacts concentrate at the lower and left side of the reconstruction. Moreover, the exponential window shows strong artifacts especially in the phase of the reconstruction as it does not fit the window form given in applications.

In our next experiment we compare the proposed tight projector against the pattern projector used in \cite{IwPrSaVi18}. We already know that both projection spaces coincide if all shifts are taken into account. Thus, we now consider the case where only every $\kappa$ shift is used for the reconstruction. Therefore let $\mathcal{T}_{LK}^{\kappa}$ be the projector onto the subspace 
$$\spa\{T_{\kappa\ell}\bar{\zb w}_k\bar{\zb w}_k^\ast T_{\kappa\ell}^\ast:\ell=0,\dots, L,\ k=0,\dots, K\},$$
and $\mathcal{P}_{LK}^{\kappa}$ the corresponding pattern projector. Note that the pattern projector does not necessarily return Hermitian matrices as required for the angular synchronization. Thus, the algorithm has to be extended to include an update step $\zb Z\leftarrow (\zb Z+\zb Z^*)/2$. A careful analysis regarding this approximation with respect to the phase retrieval problem was made by Iwen et.al. in \cite{IwPrSaVi18}.

\begin{figure}[ht]
	\center{
		\begin{tabular}{cc}
			\includegraphics[height=3.2cm]{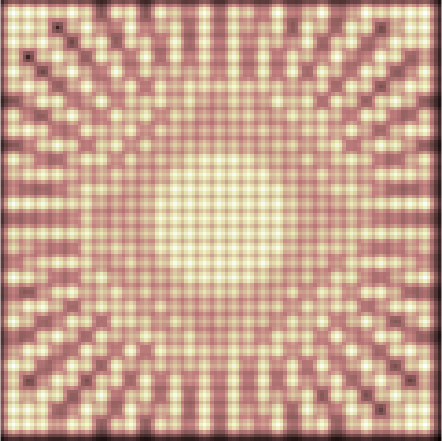} &
			\includegraphics[height=3.2cm]{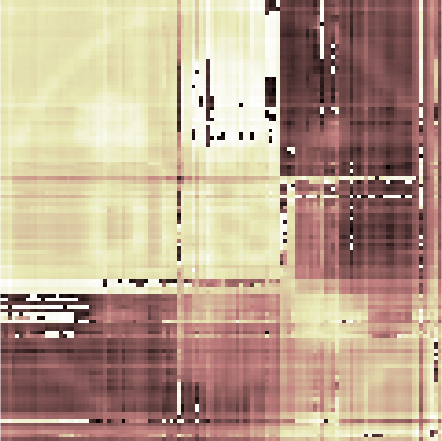} \\
			\includegraphics[height=3.2cm]{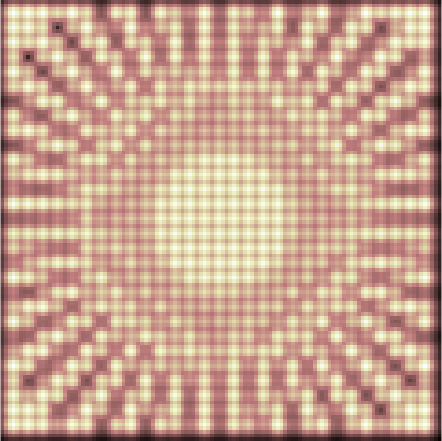} &
			\includegraphics[height=3.2cm]{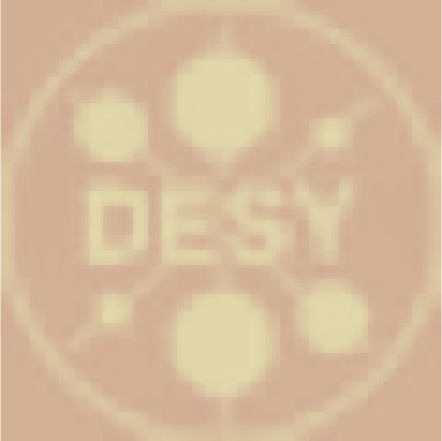}
	\end{tabular}}
	\caption{Reconstructed amplitude (left) and phase (right) using the pattern projector (top) and tight projector (bottom)}
	\label{fig:Jump4}
\end{figure}

For the experiment, we set $L$, $K$ as above, and $\kappa=4$. The results are illustrated in Figure \ref{fig:Jump4}. For both projectors a Gaussian window {\eqref{eq:GW_coef}} with $\alpha=0.99$ is used. As seen before, the reconstructed amplitude of both projectors is similar. However, the phase reconstruction is much more stable using the proposed tight projector. The pattern projector reconstruction shows strong artifacts almost dominating the original phase.

\begin{figure}[ht]
	\center{
	
			\includegraphics[height=4.5cm]{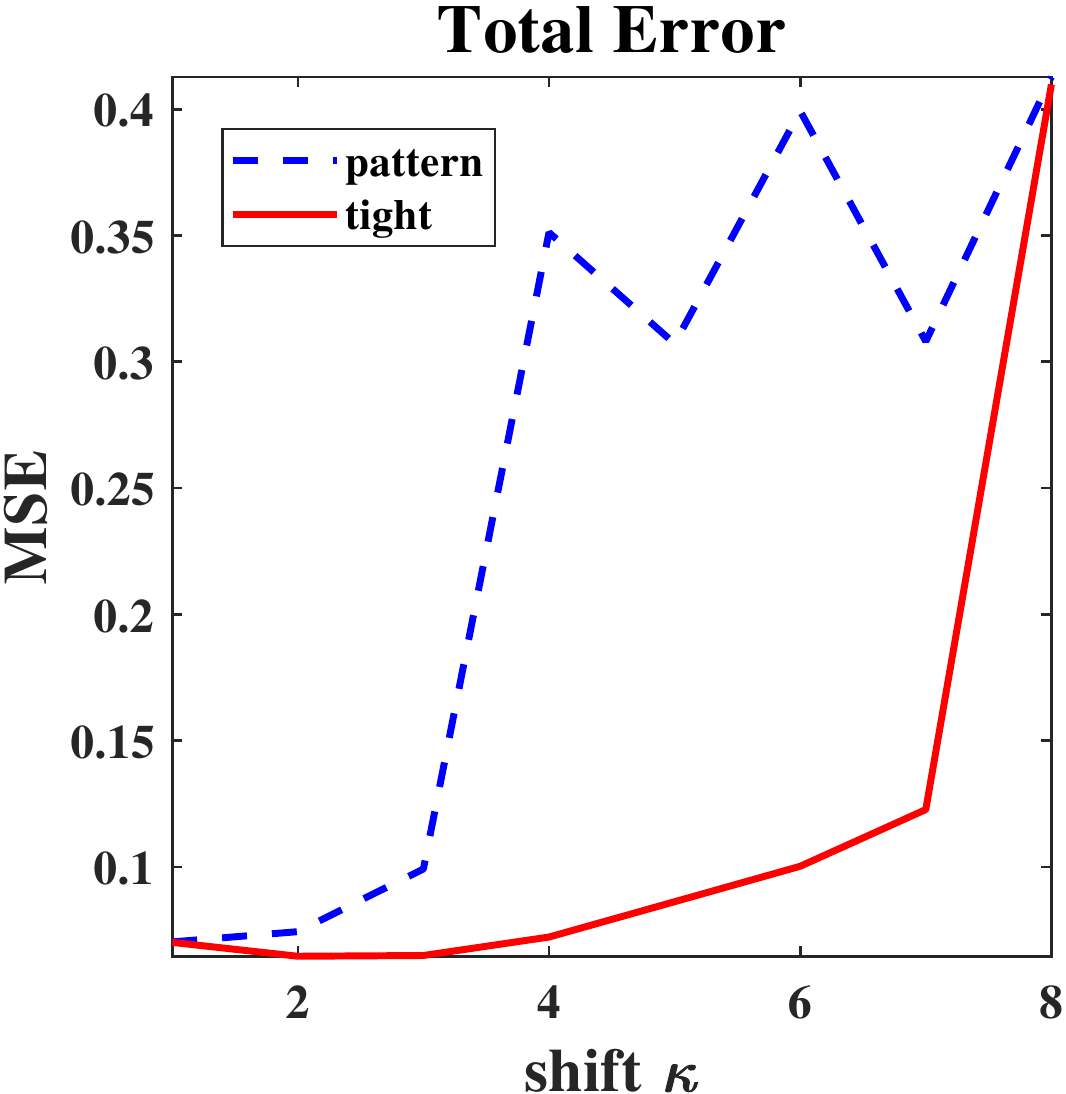} \qquad
			\includegraphics[height=4.5cm]{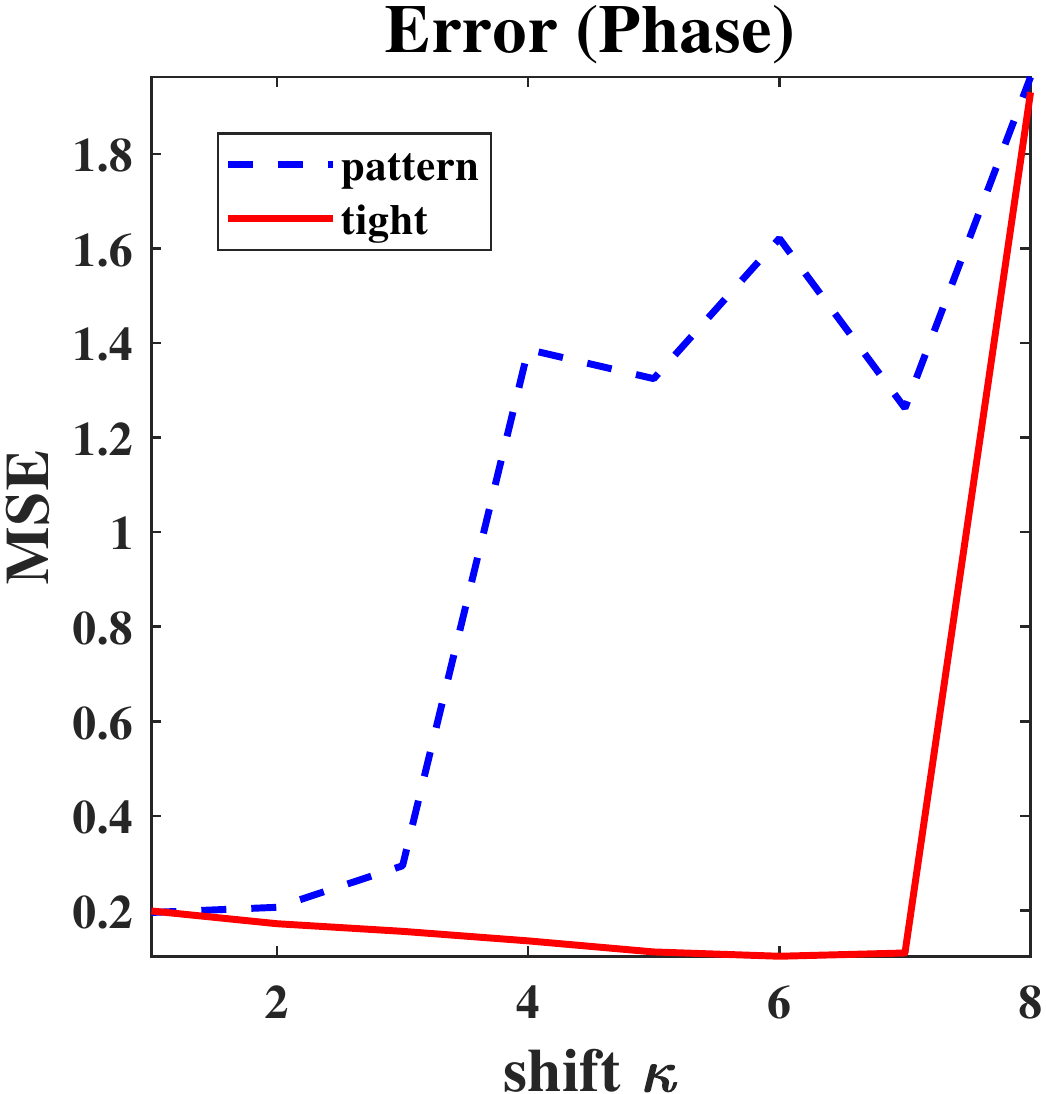} 
	
	}
	\caption{Reconstruction error using pattern and tight projector for different shift sizes $\kappa$. Left: total error, right: phase error.}
	\label{fig:projError}
\end{figure}

We verified these observations with extensive numerical experiments applying the reconstruction technique with both projectors using different window functions. {We simulated measurements of the object in Figure \ref{fig:originals} using seven different Gaussian windows. For the reconstruction we used four window functions, one exponential window \eqref{eq:2.11} and three Gaussian windows \eqref{eq:GW_coef} that differ from the windows used for simulation. Figure \ref{fig:projError} shows the mean squared error (MSE) averaged over all $28$ combinations. Here we define the MSE  as $\|\zb F-\tilde{\zb F}\|_F^2/N^2$. The original data has an amplitude range of $[0.2,0.7]$, a phase range of $[0,\pi/2]$ and a Frobenius norm of $\|\zb F\|_F=65.31$.} Clearly, the tight projector leads to a much more stable technique. As was expected, when the shift is equal to the support size of the windows (i.e., $\kappa=s=8$) both methods fail. In Figure \ref{fig:projError} the total error and the error considering only the reconstruction of the phase is shown.
The error in amplitude is not illustrated since it basically coincides for both methods. Thus, the tight operator especially stabilizes the phase reconstruction.
		
\begin{figure}[ht]
	\center{
	\includegraphics[height=4.5cm]{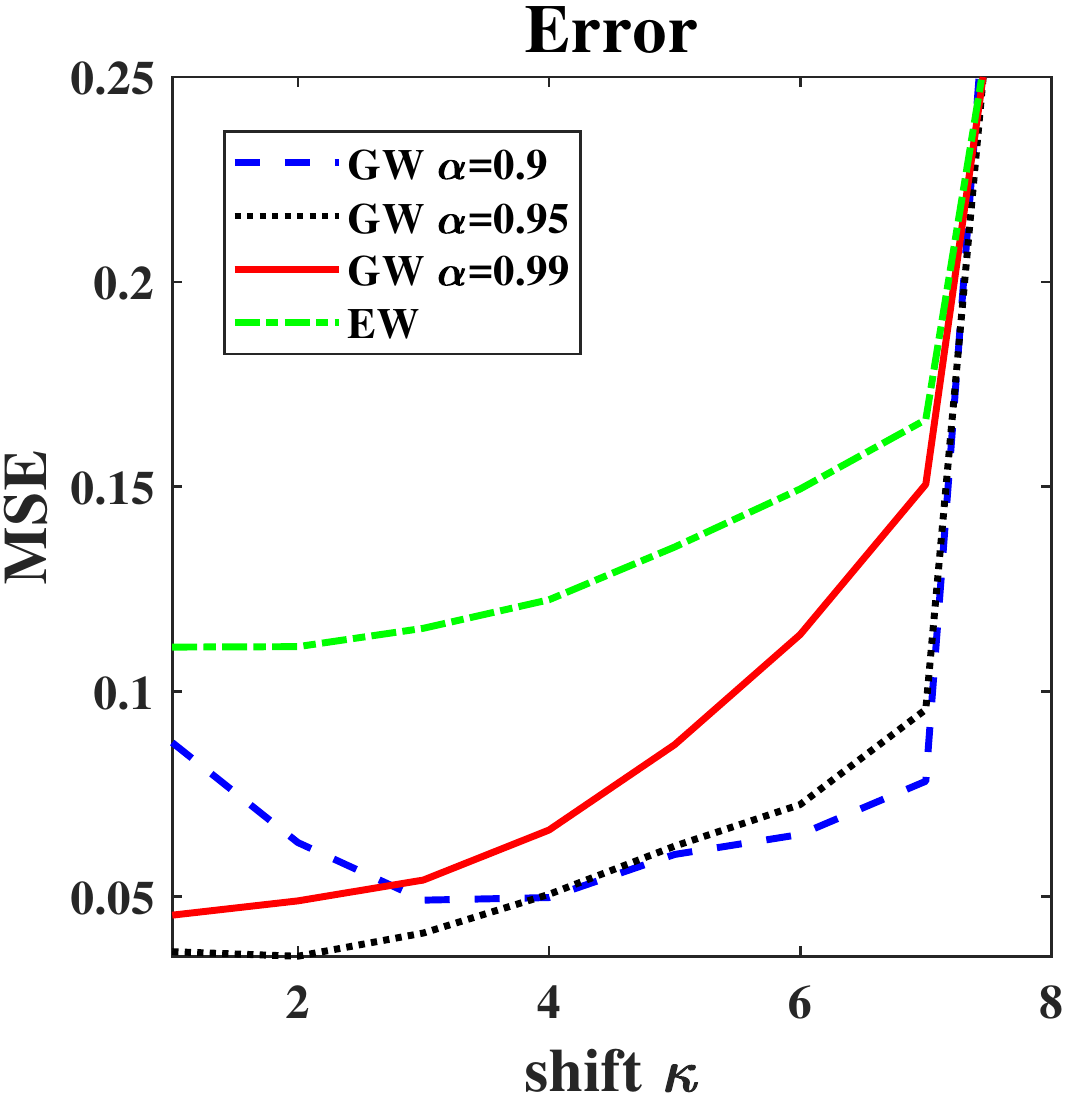}
	}
	\caption{Reconstruction error for different types of window functions using the tight projector.}
	\label{fig:windowError}
\end{figure}

Next, we compare the reconstruction error of the
tight projector for different types of window functions. The mean squared error for different shift sizes $\kappa$ is shown in Figure \ref{fig:windowError}. Here we used Gaussian windows {\eqref{eq:GW_coef}} with $\alpha=0.9,0.95,0.99$ and an exponential window {\eqref{eq:2.11}}. We observe that the reconstruction is stable up to a shift of $\kappa=s=8$ independent of the window. {The clipped Gaussian windows \eqref{eq:GW_coef} result in a smaller error since they more closely approximate the ptychographic windows used to generate the
measurements}.

Note that some window functions, such as GW with $\alpha=0.9$, even seem to perform better when not all shifts are taken into account. (Also compare the results shown in Figure \ref{fig:projError}.) This is {possibly due to the the hard cut-off of the window function for small parameters $\alpha<0.95$ which appears to contribute to Gibbs-like oscillations in the reconstructions}. 
This effect appears to diminish when the shift increases slightly.	

%As for the synchronization step the matrix has to be Hermitian we enforce this by using 
%\begin{equation}\label{eq:2.10}
%X=\frac{1}{2}\Big(\mathcal{A}_{LK}^{-1}(\zb y)+\big(\mathcal{A}_{LK}^{-1}(\zb y)\big)^\ast\Big),\quad \zb y=(y_{\ell,k}).
%\end{equation}
%for the synchronization step.
%%%%%%%%%%%%%%%%%%%%%%%%%%%%%%%%%%%%%%%%%%%%%%%%%%%%%%%%%%%%%%%%%%%%%%%%%%%%%%%%%%%%%%%%%%%
\section{\label{sec:Future}Conclusion}
%%%%%%%%%%%%%%%%%%%%%%%%%%%%%%%%%%%%%%%%%%%%%%%%%%%%%%%%%%%%%%%%%%%%%%%%%%%%%%%%%%%%%%%%%%%
We presented a direct algorithm for ptychographic reconstruction for known
windows. We have shown numerically that using a window based on the
normal distribution results in good reconstructions even if the real window
is not fully known but is assumed to be approximately Gaussian.  In addition, although the algorithm was originally designed for reconstructions based on full circulant
shifts, i.e., $L=N-1$, we have also shown that our proposed modifications result in good 
reconstructions of phase and amplitude when fewer shifts are used.
%%%%%%%%%%%%%%%%%%%%%%%%%%%%%%%%%%%%%%%%%%%%%%%%%%%%%%%%%%%%%%%%%%%%%%%%
\section*{Acknowledgements}
Mark Iwen was supported in part by NSF CCF 1615489. Rayan Saab was supported in part by NSF DMS 1517204. Nada Sissouno acknowledges support by the German Science Foundation (DFG) in the
context of the Emmy-Noether-Junior Research Group \emph{Randomized Sensing and Quantization
of Signals and Images} (KR 4512/1-1).

%%%%%%%%%%%%%%%%%%%%%%%%%%%%%%%%%%%%%%%%%%%%%%%%%%%%%%%%
\section*{References}

\bibliography{ptycho-lit}

\end{document}